%% file: m4-17.tex
\documentclass{gtart}

\input gtmonout

\volumenumber{4}
\volumename{Invariants of knots and 3-manifolds (Kyoto 2001)}
\volumeyear{2002}
\papernumber{17}  
\pagenumbers{263}{279}
\received{3 December 2001}
\revised{1 April 2002}
\accepted{22 July 2002}
\published{13 October 2002}

\usepackage{amssymb,amsmath,amscd}

\newtheorem{prob}{Problem}

\theoremstyle{definition}

\newcommand{\bp}{\begin{prob}\vskip 1ex}
\newcommand{\ep}{\end{prob}\vskip -1ex}
\newcommand{\rks}{\noindent {\bf Remarks}\qua}

\def\C{\ifmmode{\mathbb C}\else{$\mathbb C$}\fi}
\def\Z{\ifmmode{\mathbb Z}\else{$\mathbb Z$}\fi}
\def\R{\ifmmode{\mathbb R}\else{$\mathbb R$}\fi}
\def\H{\ifmmode{\mathbb H}\else{$\mathbb H$}\fi}

\newcommand{\al}{\alpha}
\newcommand{\Ga}{\Gamma}
\newcommand{\La}{\Lambda}
\newcommand{\om}{\omega}
\newcommand{\Om}{\Omega}
\newcommand{\Th}{\Theta}

\newcommand{\End}{\mathrm{End}}
\newcommand{\Sym}{\mathrm{Sym}}
\newcommand{\vol}{\mathrm{vol}}
\newcommand{\rank}{\mathrm{rank}}

\begin{document}

\title{Generalisations of Rozansky-Witten invariants} 

\authors{Justin Roberts\\Justin Sawon} 

\address{Department of Mathematics, UC San Diego\\9500 
Gilman Drive,La Jolla, CA 92093,
USA\\{\rm and}\\Department of Mathematics, SUNY at Stony Brook\\Stony Brook,
NY 11794-3651, USA}

\asciiaddress{Department of Mathematics, UC San Diego\\9500 
Gilman Drive,La Jolla, CA 92093,
USA\\and\\Department of Mathematics, SUNY at Stony Brook\\Stony Brook,
NY 11794-3651, USA}

\email{justin@euclid.ucsd.edu, sawon@math.sunysb.edu}                     
\url{www.math.ucsd.edu/\char'176justin, www.math.sunysb.edu/\char'176sawon/}

\begin{abstract}
We survey briefly the definition of the Rozansky-Witten invariants,
and review their relevance to the study of compact hyperk{\"a}hler
manifolds. We consider how various generalisations of the invariants
might prove useful for the study of non-compact hyperk{\"a}hler
manifolds, of quaternionic-K{\"a}hler manifolds, and of relations
between hyperk\"ahler manifolds and Lie algebras. The paper concludes
with a list of additional problems.
\end{abstract}

\asciiabstract{
We survey briefly the definition of the Rozansky-Witten invariants,
and review their relevance to the study of compact hyperkahler
manifolds. We consider how various generalisations of the invariants
might prove useful for the study of non-compact hyperkahler
manifolds, of quaternionic-Kahler manifolds, and of relations
between hyperkahler manifolds and Lie algebras. The paper concludes
with a list of additional problems.}

                \asciikeywords{Hyperkahler manifold, quaternion-Kahler manifold}                    

\maketitle

\section{Introduction}

In 1996 Rozansky and Witten~\cite{rw97} discovered that a {\em
hyperk\"ahler manifold} $X^{4k}$ gives rise to a {\em topological
quantum field theory} $Z_X$ in dimension $(2+1)$. In particular, to
each closed oriented $3$-manifold $M$ is associated a {\em partition
function} $Z_X(M)$, a complex-valued topological invariant of $M$,
which is defined by means of a path integral over the space of maps $M
\rightarrow X$.

This path-integral definition is not mathematically rigorous, but
Rozansky and Witten were able to give a rigorous formulation of the
(exact) {\em perturbative expansion} of $Z_X(M)$. This is an expansion
in terms of Feynman diagrams, which in this theory are simply
trivalent graphs $\Ga$:
\begin{eqnarray}
Z_X(M)=\sum_{\Ga}b_{\Ga}(X)I_{\Ga}(M).
\label{Feynman_expansion}
\end{eqnarray}
The expansion decouples the partition function into quantities
$b_{\Ga}(X)$ depending only on $X$, and quantities $I_{\Ga}(M)$
depending only on $M$. Habegger and Thompson~\cite{ht99} showed that
when $b_1(M)\geq 1$, the numbers $I_{\Ga}(M)$ are essentially
coefficients of the {\em Le-Murakami-Ohtsuki invariant}~\cite{lmo98}
of $M$; this is believed to hold also when $b_1(M)$ is zero. Therefore
the new ingredient in Rozansky and Witten's construction is the {\em
weight system} $b_{\Ga}(X)$, a complex-valued function on trivalent
graphs, constructed from the geometry of $X$.

In this paper we will investigate the weight systems primarily {\em as
geometers}, with a view to using them to study the manifold $X$.
Although it is also very interesting to study the invariants from the
point of view of Vassiliev theory \cite{rbw01}, or of TQFT
\cite{rsw01}, we will not pursue these ideas here.

Let us review the definition of the Rozansky-Witten weight systems,
assuming for the moment that $X$ is {\em compact}. Recall that a
{\em hyperk{\"a}hler manifold} $X^{4k}$ is a Riemannian manifold admitting
three orthogonal parallel complex structures $I,J,K$ which satisfy
the quaternion identities. Its holonomy group is contained in
$\mathrm{Sp}(k)$, and therefore the frame bundle can be reduced to a
principal $\mathrm{Sp}(k)$-bundle. The complexification of the tangent
bundle therefore splits as a tensor product
\begin{eqnarray}
TX\otimes_{\R}\C =E\otimes_{\C} H
\label{EtimesH}
\end{eqnarray}
where $E$ is the bundle induced by the standard representation of
$\mathrm{Sp}(k)$ on $\H^k$ ($\H$ denoting the quaternions) and $H$ is a
trivialisable complex rank-two bundle. (This notation is consistent
with Salamon's~\cite{salamon99} but differs from that of Rozansky and
Witten~\cite{rw97}.) 

The Riemann curvature tensor reduces to a section $\Psi$ of
$\mathrm{Sym}^4E$, and the metric induces skew-symmetric two-forms on
$E$ and $H$. Suppose we are given a trivalent graph $\Ga$ with $2k$
vertices. We place a copy of $\Psi$ at each vertex, labelling the
outgoing edges by the indices of $\Psi$ (with one extra index at each
vertex). Then we tensor these copies of $\Psi$ together, contract
indices along edges using the skew two-form on $E$, and project the
result from $\otimes^{2k}E$ to $\La^{2k}E$. The latter bundle has a
trivialisation given by the $k$th power of the skew two-form, and
hence we obtain a well-defined scalar field which we can integrate
over $X$. Up to normalization (including orientation conventions to
handle the signs) this gives $b_{\Ga}(X)$, an invariant of the
hyperk{\"a}hler structure on $X$ which we will refer to as a {\em
Rozansky-Witten invariant\/} of $X$.

Rozansky and Witten showed that $b_{\Ga}(X)$ depends on the graph
$\Ga$ only up to the {\em anti-symmetry} and {\em IHX relations}, and
hence determines a weight system on the {\em graph homology} space of
rational linear combinations of equivalence classes of graphs --- that
is, the space of Jacobi diagrams typically denoted elsewhere in these
proceedings by ${\mathcal A}(\emptyset)$. In particular, this implies
that $Z_X(M)$ really is a topological invariant of $M$, for the terms
$I_{\Ga}(M)$ actually depend on a choice of metric on $M$ and vary as
the metric is deformed, but the variations in the
sum~(\ref{Feynman_expansion}) cancel out.

If we choose a complex structure on $X$ which is compatible with the
hyperk{\"a}hler metric, then the resulting complex manifold has a
holomorphic symplectic form $\omega\in\mathrm{H}^0(X,\La^2T^*)$, where
$T^*$ is the holomorphic cotangent
bundle. Kontsevich~\cite{kontsevich99} gave a reinterpretation of the
weights $b_{\Ga}(X)$ using characteristic classes of the resulting
complex-valued symplectic foliation. At the same time
Kapranov~\cite{kapranov99} recognised that from this point of view,
the curvature tensor can be viewed as a Dolbeault representative of
the {\em Atiyah class\/} $\alpha_T$ of the holomorphic tangent
bundle. The Atiyah class $\alpha_F$ of a holomorphic bundle $F$ is the
extension class of the short exact sequence $$0\rightarrow T^*\otimes
F\rightarrow
\mathrm{J}^1(F)\rightarrow F\rightarrow 0$$ where $\mathrm{J}^1(F)$ is
the first jet bundle of $F$ (see~\cite{atiyah57}). Thus the Atiyah
class $\alpha_F$ is an element of $\mathrm{Ext}^1(F,T^*\otimes
F)=\mathrm{H}^1(X,T^*\otimes\End (F))$; it vanishes if and only if the
sequence splits, and hence $\alpha_F$ is the obstruction to the
existence of a global holomorphic connection on $F$. We can define the
weights $b_{\Ga}(X)$ using the Atiyah class instead of the curvature
tensor, with the advantage that we no longer need to know the metric
on $X$ (usually we know only of the {\em existence\/} of a
hyperk{\"a}hler metric, by Yau's Theorem). This approach required us
to choose a complex structure on $X$, but we already know from the
original definition that $b_{\Ga}(X)$ is independent of this choice,
which amounts to a choice of holomorphic trivialisation of the bundle
$H$, and gives an identification of $E$ with the holomorphic tangent
bundle $T$. Consequences of this approach are that Rozansky-Witten
invariants are invariant under {\em deformation} of the hyperk\"ahler
metric on $X$, and in fact can be defined more generally for {\em
holomorphic symplectic manifolds\/}, which are hyperk\"ahler if and
only if they are K\"ahler.

The first extensive calculations of Rozansky-Witten invariants were
carried out by the second author in~\cite{sawon99}. Let $\Th^k$
be the (disconnected) graph given by $k$ copies of the {\em theta
graph\/} $\Th$. If $X$ is an {\em irreducible\/} hyperk{\"a}hler
manifold (i.e.\ no covering space of $X$ splits into a product) then
$b_{\Th^k}(X)$ can be written
\begin{eqnarray}
\frac{2^kk^k}{k!}\frac{(\int_Xc_2[\omega\bar{\omega}]^{k-1})^k}{(\int_X[\omega\bar{\omega}]^k)^{k-1}}=\frac{k!}{(4\pi^2k)^k}\frac{\|R\|^{2k}}{(\vol(
X)
)^{k-1}}
\label{Thetak}
\end{eqnarray}
where $c_2$ is the second Chern class of $T$, $[\om\bar{\om}]$ is the
cohomology class represented by $\om\bar{\om}$, $\|R\|$ is the
$L^2$-norm of the curvature, and $\vol (X)$ is the volume. (This is
equation (10) in~\cite{hs01}.) Examples of irreducible hyperk{\"a}hler
manifolds of dimension $4k$ are the {\em Hilbert scheme} of $k$ points
on a K3 surface and the {generalized Kummer varieties}
(see Beauville~\cite{beauville83}).

The {\em Chern numbers} of $X$ can also be expressed (see
\cite{sawon99}) as Rozansky-Witten invariants $b_{\Ga}(X)$, where we
choose $\Ga$ to be the {\em closure of a disjoint union of wheel
graphs\/}. Not all graphs can be expressed as a linear combination of
such closures, and indeed Rozansky-Witten invariants are strictly more
general than characteristic numbers (see~\cite{sawon00i}, page
360). However, it is a surprising fact that $\Th^k$ can be expressed
in such a way: this is a simple corollary of the Wheeling Theorem of
Bar-Natan, Le, and Thurston~\cite{thurston00}; a more direct proof was
recently given by Britze and Nieper~\cite{bn01}. We find that
\begin{eqnarray}
b_{\Th^k}(X)=48^kk!\hat{A}^{1/2}[X]
\label{Ahat}
\end{eqnarray}
where $\hat{A}^{1/2}$ is the square root of the
$\hat{A}$-polynomial. Combining this with the
expression~(\ref{Thetak}) for $b_{\Th^k}(X)$ gives us a formula (the
main result of~\cite{hs01}) for the $L^2$-norm of the curvature in
terms of characteristic numbers and the volume of $X$. This formula
was an important ingredient in Huybrechts' proof \cite{huybrechts01}
of a finiteness result for the number of diffeomorphism types of
compact hyperk{\"a}hler manifolds.

There is a useful generalisation \cite{kapranov99, sawon00i,
thompson99} of the basic construction. A {\em holomorphic bundle} (or
coherent sheaf) $F$ on $X$ allows us to define numbers $b_{\Ga}(X,F)$
associated to any trivalent graph $\Ga$ which has a preferred oriented
circle. This is the procedure by which one obtains genuine Vassiliev
weight systems, defined on the usual algebra ${\mathcal A}$ of Jacobi
diagrams.

%%%%%%%%%%%%%%%%%%%%%%%%%%%%%%%%%%%%%%%%%%%%%%%%%%%%%%%%%%%%%%%%%%%%%
%%%%%%%%%%%%%%%%%%%%%%%%%%%%%%%%%%%%%%%%%%%%%%%%%%%%%%%%%%%%%%%%%%%%%

\section{Non-compact hyperk{\"a}hler manifolds}

There are several families of {\em non-compact} hyperk{\"a}hler
manifolds $X$ whose Rozansky-Witten invariants are interesting for
physical reasons. 

{\bf (a) Monopole spaces}\qua The story of the invariants actually
began with a dimensionally reduced version of Seiberg-Witten theory
with gauge group $\mathrm{SU}(2)$. There is a topologically twisted
version of this theory which has as its partition function the {\em
Casson invariant}, by which we mean the generalisation due to Walker
and Lescop~\cite{lescop96} which works for arbitrary 3-manifolds, and
not just for homology spheres. In the low energy limit, the theory
becomes a sigma-model with target space the {\em Atiyah-Hitchin
manifold} $X_{\mathrm{AH}}$, which is the reduced moduli space of
{\em two-monopoles} on $\R^3$, and is a smooth non-compact hyperk{\"a}hler
manifold of dimension four. Rozansky and Witten observed that in fact
their model was still defined with $X_{\mathrm{AH}}$ replaced by any
(compact or asymptotically flat) hyperk\"ahler manifold.

There should be versions of the Casson invariant for gauge groups
other than $\mathrm{SU}(2)$. For a general group $G$ we should replace
$X_{\mathrm{AH}}$ by the moduli space of vacua of the 3-dimensional
SUSY gauge theory with gauge group $G$; for example, when $G$ is
$\mathrm{SU}(N)$ we should choose the target space $X$ to be the
reduced moduli space of {\em $N$-monopoles}. The resulting partition
function would be a finite-type invariant of Ohtsuki order $3(N-1)$,
and it would be particularly interesting to compute it in the case of
$SU(3)$ (perhaps by a surgery formula - see
Thompson~\cite{thompson99}, for example), and to investigate whether
there is any connection with the generalized Casson invariant of Boden
and Herald \cite{bh}.

{\bf (b) Hilbert schemes}\qua The (reduced) Hilbert scheme of points on
$\C^2$ are hyperk\"ahler manifolds whose Rozansky-Witten theory is
supposed to arise when reducing $M$-theory on manifolds with $G_2$
holonomy to associative submanifolds. Consequently, it would be useful
to calculate $b_{\Ga}(X)$ for these manifolds.

{\bf (c) Gravitational instantons}\qua (Kronheimer's {\em ALE
spaces}~\cite{kronheimer89})\qua These spaces are resolutions of quotient
singularities $\C^2/G$, where $G$ is a finite subgroup of
$\mathrm{SU}(2)$. The McKay correspondence says that such subgroups
correspond to Dynkin diagrams of type A, D, and E. For example, the
cyclic group of order $N$ corresponds to A$_{N-1}$. Kronheimer gave a
construction of these spaces as hyperk{\"a}hler
quotients.

When $X$ is non-compact we encounter a couple of basic problems. Most
obviously, the integral in the definition of $b_{\Ga}(X)$ may not
converge; but even if it does, we won't necessarily get a weight
system on {\em graph homology\/}, as it may be necessary to add a
boundary term to the IHX relation. 

One can always work with weight systems taking graph homology to the
{\em Dolbeault cohomology} of the manifold, but this might itself be
very complicated, and therefore useless. So let us consider some
possible approaches to the convergence problem for the quantities
$b_{\Ga}(X)$. The first is a direct differential geometric approach:
there are examples of non-compact hyperk{\"a}hler manifolds for which
the metric is {\em explicitly known}, unlike in the compact
case. However, the calculations involved are not simple; see Sethi,
Stern, and Zaslow~\cite{ssz95} for this direct calculation when $X$ is
the Atiyah-Hitchin manifold. A second approach is to try to mimic
Kapranov's approach, but using {\em compactly supported} (or perhaps
$L^2$) cohomology. A third approach is to {\em compactify} the
manifold; we won't get a compact {\em hyperk{\"a}hler} manifold but we
can at least try to define Rozansky-Witten invariants of the manifold
so obtained.

This third approach has been investigated by Goto~\cite{goto01}. He
generalised holomorphic symplectic manifolds by defining {\em log
symplectic manifolds\/} as (compact) complex manifolds having a
meromorphic symplectic form with logarithmic poles along a given
divisor $D$ (which must be the anti-canonical divisor). Under certain
hypotheses, Rozansky-Witten invariants can be defined for this class
of manifolds. If $X$ is non-compact, we can try to compactify by
adding a divisor at infinity. In the case of the reduced monopole
moduli spaces this results in a log symplectic manifold. We don't yet
know whether the (log symplectic) Rozansky-Witten invariants of this
compactification agree with the Rozansky-Witten invariants of the
non-compact space, as defined in terms of the metric, but this
approach seems highly promising.

The second approach can be illustrated in the case of gravitational
instanton spaces. First note that in four dimensions there is
essentially just one Rozansky-Witten invariant
$$b_{\Th}(X)=\frac{1}{(2\pi)^2}\int_X \mathrm{Tr}(R^2)$$ which can be
calculated using the Gauss-Bonnet formula. We find
$$b_{\Th}(X)=2\chi(X)+B(X)$$ where $\chi$ is the Euler characteristic
and $B$ is a boundary term. For example, when $X$ is $\C^2$ we have
$b_{\Th}(\C^2)=0$ and $\chi(\C^2)=1$, and therefore the boundary term
$B(\C^2)$ must equal $-2$. Now gravitational instanton spaces are {\em
asymptotically locally Euclidean\/} (ALE), which means their metric is
asymptotic to the metric on $\C^2/G$ at infinity. In particular, this
means they have boundary term $$B(X)=\frac{B(\C^2)}{|G|},$$ where
$|G|$ is the order of the group $G$, which acts freely on $\C^2$ at
infinity. The homology group $\mathrm{H}_2(X,\Z)$ is generated by the
irreducible curves $C_i$ which make up the exceptional divisor $E$,
and these correspond to the nodes of the Dynkin diagram; apart from
$\mathrm{H}_0(X,\Z)\cong\Z$, all other homology groups are
zero. Therefore $$\chi(X)=1+\rank(G).$$ Combining the above, we get
$$b_{\Th}(X)=2+2\,\rank(G)-\frac{2}{|G|}.$$ In particular, for the
A-series we get $2N-\frac{2}{N}$. This happens to be twice the value
of the Casimir in the fundamental representation of $SU(N)$, but 
we don't know a general Lie-theoretic interpretation of the formula. 

Now recall that the Atiyah class $\al_T$ is the obstruction to
the existence of a global holomorphic connection. If we choose a cover
$\{U_i\}$ of $X$ and local holomorphic connections $\nabla_i$ on each
open set $U_i$, then the differences
$$\nabla_i-\nabla_j\in\mathrm{H}^0(U_i\cap U_j,T^*\otimes{\End}(T))$$
form a 1-cocycle representing $\al_T$ in {\v C}ech cohomology. Using
the ALE property again, we can choose a holomorphic (in fact flat)
connection on the complement $U_0$ of a compact neighbourhood of the
exceptional divisor $E$. This gives an Atiyah class which is {\em
compactly supported} in a neighbourhood of $E$, and this Atiyah class
can be used to calculate $b_{\Th}(X)$. In the case when $X$ is the
A-series of gravitational instantons, this calculation has been
carried out by the second author and the answer agrees with the one
given above.

In general, we don't expect there to be a compactly supported
representative of the Atiyah class on an arbitrary non-compact
hyperk{\"a}hler manifold $X$. However, if $X$ is asymptotically flat
then there ought to be a representative with some kind of nice
asymptotic behaviour (perhaps lying in $L^2$-cohomology, for
instance). We are not yet aware of any such examples.

\section{Quaternionic-K{\"a}hler manifolds}

Hyperk{\"a}hler manifolds can be regarded as particular examples of
{\em quaternionic-K{\"a}hler\/} manifolds. Whereas the former have
holonomy contained in $\mathrm{Sp}(k)$, quaternionic-K{\"a}hler
manifolds are characterized by having holonomy contained in
$$\mathrm{Sp}(k).\mathrm{Sp}(1)=\mathrm{Sp}(k)\times\mathrm{Sp}(1)/{\Z_2}$$ 
where $\Z_2$ is generated by $(-1,-1)$. The scalar curvature $s$ of
a quaternionic-K{\"a}hler manifold $X$ vanishes if and only if it is 
locally hyperk{\"a}hler (in particular, if $X$ is simply connected
then $s$ vanishes if and only if $X$ is globally
hyperk{\"a}hler). However, if $s$ is non-zero then $X$ is not
in general K{\"a}hler and does not have a holomorphic symplectic form
(indeed it may not even admit an almost complex structure). Thus at
first sight it may appear impossible to define anything resembling
Rozansky-Witten invariants for quaternionic-K{\"a}hler manifolds. Let
us assume once again that all our manifolds are compact; our results
on quaternionic-K{\"a}hler manifolds will be drawn from Salamon's
survey article~\cite{salamon99}.

For a quaternionic-K{\"a}hler manifold $X$ there is an obstruction
$\epsilon\in\mathrm{H}^2(X,\Z_2)$ to the lifting of the structure
group $\mathrm{Sp}(k).\mathrm{Sp}(1)$ of $X$ to
$\mathrm{Sp}(k)\times\mathrm{Sp}(1)$. If $\epsilon$ vanishes then the
decomposition~(\ref{EtimesH}) of the complexified tangent bundle is
still valid, but with $H$ now a {\em non-trivial\/} rank-two complex
vector bundle. Indeed $H$ is the bundle induced by the standard
representation of $\mathrm{Sp}(1)$ on $\H$, just as $E$ is induced by
the standard representation of $\mathrm{Sp}(k)$ on $\H^k$. More
generally, over open sets where $\epsilon$ vanishes, there is a {\em
local} decomposition of this form; for simplicity we assume $\epsilon$
vanishes globally.

Let us return to Rozansky and Witten's original definition of the
invariants using $\Psi$, the section of $\Sym^4E$. In the
quaternionic-K{\"a}hler case the curvature tensor of $X$ equals
$R_Q+s\rho_1$ where $R_Q$ is a section of $\Sym^4E$ and $\rho_1$ is a
certain invariant element (see~\cite{salamon99}, Corollary 3.4). We
can then proceed as before using $R_Q$ instead of $\Psi$, and thus
define Rozansky-Witten invariants $b_{\Ga}(X)$ for
quaternionic-K{\"a}hler manifolds $X$.

There ought to be some interesting consequences to this
approach. Since we are using $R_Q$ and not the full curvature tensor,
we don't expect to get characteristic numbers, as in the
hyperk{\"a}hler case. However, if the scalar curvature $s$ is of fixed
sign, we might generate {\em inequalities\/} between these
Rozansky-Witten invariants and characteristic numbers. For example, it
is conjectured that all positive (i.e.\ $s>0$) quaternionic-K{\"a}hler
manifolds are symmetric (known as Wolf spaces), and this has been 
proved in dimensions $4$, $8$, and $12$. In dimension $12$, bounds on
the $\hat{A}$-genus play an important role in the proof; it is
tempting to believe that equations like~(\ref{Ahat}), suitably
generalised to the quaternionic-K{\"a}hler case, would be useful in
tackling the conjecture in higher dimensions.

There is an alternative approach to defining invariants if we wish to
stay closer to authentic characteristic numbers. We first recall some
facts about characteristic classes of quaternionic-K{\"a}hler
manifolds (see~\cite{salamon99}, Section 8). Denote by $u$ minus the 
second Chern class $-c_2(H)$ of the bundle $H$. Then $u$ and the Chern
classes of $E$ are well-defined elements of the rational cohomology
ring of $X$ (this is true even when $\epsilon$ is non-zero). Moreover,
$4u$ is actually integral. From
$$\mathrm{ch}(TX\otimes_{\R}\C)=\mathrm{ch}(E)\mathrm{ch}(H)$$
we see that the first Pontrjagin class of $X$ is
$$p_1=2(ku-c_2(E)).$$

If we look at the left hand side of equation~(\ref{Thetak}) we see
that we only need $c_2$ and the $4$-form $\om\bar{\om}$ to define
$b_{\Th^k}(X)$. Recall that for a hyperk{\"a}hler manifold the odd
Chern classes vanish and the even Chern classes are equivalent to the
Pontrjagin classes, so we can replace $c_2$ by $p_1$. On a
quaternionic-K{\"a}hler manifold $X$ there is a natural $4$-form
$\Om$ which is parallel and hence closed. One way to define $\Om$ is
to take local almost-complex structures $\{I,J,K\}$ behaving like the
quaternions ($IJ=K$), and corresponding $2$-forms defined by
$\om_I(v,w)=g(Iv,w)$, etc. Then
$$\Om=\om_I\wedge\om_I+\om_J\wedge\om_J+\om_K\wedge\om_K$$
is well-defined up to overall scale, and parallel (even though the
local almost-complex structures and $2$-forms are not). Thus in the 
quaternionic-K{\"a}hler case we can replace $\om\bar{\om}$ in
equation~(\ref{Thetak}) by $\Om$, and hence define
$$b_{\Th^k}(X)=\frac{2^kk^k}{k!}\frac{(\int_Xp_1[\Om]^{k-1})^k}{(\int_X[\Om]^k)^{k-1}}.$$
The right hand side is invariant under rescaling of $\Om$. In fact
there is a normalised version of $\Om$ (involving the scalar curvature
$s$) which represents the cohomology class $u$, so the definition
above is equivalent to
$$b_{\Th^k}(X)=\frac{2^kk^k}{k!}\frac{(\int_Xp_1u^{k-1})^k}{(\int_Xu^k)^{k-1}}.$$  
The denominator is known as the quaternionic volume, and is
positive.

Clearly this expression is just one of a family of expressions
involving integrals of Pontrjagin classes and the cohomology class
$u$. In the hyperk{\"a}hler case (using Chern numbers and
$[\om\bar{\om}]$) all of these expressions are Rozansky-Witten
invariants, though there is no reason to expect that the converse is
true. So it may be that we don't get a full generalisation of
Rozansky-Witten invariants to quaternionic-K{\"a}hler manifolds in
this way. Perhaps more importantly, we no longer have the machinery of
graph homology which enabled us to prove relations such as
equation~(\ref{Ahat}).

%%%%%%%%%%%%%%%%%%%%%%%%%%%%%%%%%%%%%%%%%%%%%%%%%%%%%%%%%%%%%%%%%%%%%
%%%%%%%%%%%%%%%%%%%%%%%%%%%%%%%%%%%%%%%%%%%%%%%%%%%%%%%%%%%%%%%%%%%%%

\section{Hyperk\"ahler manifolds and Lie algebras}

%%%%%%%%%%%%%%%%%%%%%%%%%%%%%%%%%%%%%%%%%%%%%%%%%%%%%%%%%%%%%%%%%%%%%
%%%%%%%%%%%%%%%%%%%%%%%%%%%%%%%%%%%%%%%%%%%%%%%%%%%%%%%%%%%%%%%%%%%%%

The primary examples of Vassiliev weight systems are those coming from
complex semisimple Lie algebras and superalgebras. It was once
conjectured that in fact {\em all} weight systems were of this form,
but this was shown to be false by Vogel \cite{Vogel}. 

From the point of view of Vassiliev theory, then, the most obvious
question is whether the Rozansky-Witten weight systems are really {\em
new}, that is, lying outside the span of the Lie algebraic weight
systems. To resolve this issue it would probably be necessary to
understand how the RW weight systems behave under the action of
Vogel's algebra $\Lambda$.

A more general problem is to understand the common properties of
hyperk\"ahler manifolds and Lie algebras. There is in fact an {\em
algebraic} way to unite them, described in \cite{rbw01}, where it is
shown that (roughly) the derived category of coherent sheaves on
$X^{4n}$ is the category of modules over the holomorphic tangent sheaf
$T$, which is a Lie algebra in this category. Thus, a hyperk\"ahler
manifold {\em can} be thought of as giving rise to a Lie algebra --
but not in the usual category of vector spaces.

Studying the analogy between the derived category and a representation
category seems fruitful. There are similiarities: the derived category
possesses a Duflo isomorphism, namely ``theorem of the M.~Kontsevich on
complex manifold'' \cite{kontsdefquant}. But there are also basic
differences: the derived category is not semisimple. It would be
interesting to study whether there are analogues of Adams
operations, character formulae, Harish-Chandra isomorphisms, etc. in
the derived category; whether there are classes of sheaves
(exceptional, stable, etc.) or bundles which behave like irreducible
objects; and whether there is a meaningful interpretation of the
``quantum dimensions'' of sheaves arising from the Rozansky-Witten
partition function of the unknot.

In the converse direction, is it possible to realise the Lie algebraic
weight systems as the Rozansky-Witten invariants of suitable
hyperk\"ahler manifolds? These would presumably have to be {\em
infinite-dimensional} manifolds, because a Lie algebra defines
weight systems in all degrees, whereas a manifold $X^{4k}$ defines
weight systems only in degree less than or equal to $k$. A potentially
useful observation is that the based loop space $\Omega G^{\C}$ of a
{\em complex} semisimple Lie group is a hyperk\"ahler manifold; the
idea to look at the loop space (or perhaps the classifying space $BG$)
is also suggested by Kontsevich's description \cite{kontsevich99} of a
``shifted'' Lie algebra as a formal symplectic manifold. The formal
shift in grading is closely related to the shift used in forming the
bar complex to obtain a model of the cohomology of the loop space.

Returning again to gravitational instantons, we have seen a tiny piece
of numerical evidence suggesting that the invariants $b_{\Ga}(X)$
corresponding to these spaces should be related to the weight systems
of the corresponding Lie algebras (i.e.\ of types A, D, and E). Of
course the ALE spaces are four-dimensional and so there is only one
invariant $b_{\Th}(X)$; to obtain invariants in all degrees we would
need to consider also the Hilbert schemes of points on $X$, or more
generally moduli spaces of instantons on $X$. These spaces can also be
obtained by the hyperk{\"a}hler quotient construction: see Kronheimer
and Nakajima~\cite{kn90}.

A slightly different way to look for connections between hyperk\"ahler
manifolds and Lie algebras is to work in terms of the associated
TQFTs, rather than just weight systems. Each TQFT is determined by an
underlying {\em ribbon category\/} $\cal C$. In Witten's Chern-Simons
theory \cite{wittencs}, $\cal C$ is the category of representations of
a quantum group at a root of unity (see Reshetikhin and Turaev
\cite{RT}), whereas in Rozansky-Witten theory we find \cite{rsw01}
that $\cal C$ is essentially the derived category of coherent sheaves
on $X$, with a non-standard ribbon category structure. Might there
then be a manifold whose derived category is one of these
representation categories? This seems somewhat outrageous, because one
would not expect the derived categories to be semisimple, or even
abelian. One could at least hope for some extension of the work of
Kapranov and Vasserot \cite{kapvass}, who described the derived
category of the ALE spaces in purely algebraic terms. Note also that
Freed, Hopkins and Teleman \cite{fht} have shown that the
representation ring (Verlinde algebra) of the quantum group appears
geometrically as a twisted $K$-group; it is tempting to wonder whether
the underlying representation category might not arise as some kind of
underlying twisted derived category of sheaves (such things have been
defined by Caldararu \cite{cald}). 

A totally different approach to connecting the two worlds would be to
define a theory of {\em equivariant} Rozansky-Witten invariants, in
which the $G$-equivariant invariants of a point give rise to Lie
algebra weight systems. Quite what properties such a theory should
have is unclear: it should perhaps be related more to the
hyperk\"ahler quotient construction than to the usual homotopy
quotient underlying equivariant cohomology.

%%%%%%%%%%%%%%%%%%%%%%%%%%%%%%%%%%%%%%%%%%%%%%%%%%%%%%%%%%%%%%%%%%%%%
%%%%%%%%%%%%%%%%%%%%%%%%%%%%%%%%%%%%%%%%%%%%%%%%%%%%%%%%%%%%%%%%%%%%%

\section{Further problems}

%%%%%%%%%%%%%%%%%%%%%%%%%%%%%%%%%%%%%%%%%%%%%%%%%%%%%%%%%%%%%%%%%%%%%
%%%%%%%%%%%%%%%%%%%%%%%%%%%%%%%%%%%%%%%%%%%%%%%%%%%%%%%%%%%%%%%%%%%%%

\bp Derive new constraints of the form {\rm (}\ref{Ahat}{\rm )} for
hyperk\"ahler manifolds. \ep

\rks We hope that the existence of the TQFT structure, which allows
the computation of partition functions in many different ways, will
yield new identities involving the characteristic classes of
hyperk\"ahler manifolds. It would be particularly useful to give a
formula for the norm of the curvature of a general holomorphic bundle,
bounding the Chern numbers of such bundles on $X$.

%%%%%%%%%%%%%%%%%%%%%%%%%%%%%%%%%%%%%%%%%%%%%%%%%%%%%%%%%%%%%%%%%%%%%

\bp Are the invariants $b_{\Ga}(X)$ rational, or even integral? \ep

\rks Current evidence suggests that, correctly normalised, the weight
systems for compact hyperk\"ahler manifolds take {\em integral} values
on trivalent diagrams, and are at least rational for non-compact
manifolds. Might the integrality be a consequence of an alternative,
enumerative-geometric ``counting'' definition of the invariants?

%%%%%%%%%%%%%%%%%%%%%%%%%%%%%%%%%%%%%%%%%%%%%%%%%%%%%%%%%%%%%%%%%%%%%

\bp What is the meaning of the $\hat A^{\frac12}$-genus of a manifold?
\ep

\rks This genus appears in the identity (\ref{Ahat}). The $\hat
A^{\frac12}$-class appears also in Kontsevich's Duflo isomorphism for
complex manifolds, and in physics, in the study of $D$-brane charge
\cite{dbrane}. One would assume that these occurrences are somehow
related. Is there some natural class of manifolds for which the genus
is integral? In other words, might it be the index of an elliptic
operator which can only be defined under certain geometrical
conditions? (It is {\em not} actually integral for hyperk\"ahler
manifolds --- see \cite{sawon99}.)

%%%%%%%%%%%%%%%%%%%%%%%%%%%%%%%%%%%%%%%%%%%%%%%%%%%%%%%%%%%%%%%%%%%%%

\bp What is the relation between RW invariants and Chern numbers? \ep

\rks We can consider various classes of invariants, each a
generalisation of the previous one: linear combinations of honest
Chern numbers; generalised Chern numbers involving powers of the
symplectic form and its conjugate; RW invariants; ``big Chern
classes'' \cite{kapranov99}. The linear combinations of Chern numbers
are contained in rational functions of Chern numbers, which are
determined by the Chern classes. Which of these inclusions is strict?
Are there other interesting ones?

%%%%%%%%%%%%%%%%%%%%%%%%%%%%%%%%%%%%%%%%%%%%%%%%%%%%%%%%%%%%%%%%%%%%%

\bp Explain the ``characteristic numbers'' of non-compact
examples. \ep

\rks The partition functions $Z_X(S^1 \times S^2)$ and $Z_X(S^1 \times
S^1 \times S^1)$ are, for any TQFT, the dimensions of certain vector
spaces, and hence integral. In RW theory on a compact manifold $X$,
they are the {\em Todd genus} and {\em Euler characteristic} of $X$,
respectively. For non-compact $X$ they seem to be rational: for
example $Z_{X_{\mathrm{AH}}}(S^1 \times S^2)=-\frac{1}{12}$ according
to Rozansky and Witten (at least up to sign, about which there are
some difficulties). Might such numbers in fact be the regularised
dimensions of infinite-dimensional vector spaces?

%%%%%%%%%%%%%%%%%%%%%%%%%%%%%%%%%%%%%%%%%%%%%%%%%%%%%%%%%%%%%%%%%%%%%

\bp Sen's conjecture on the cohomology of monopole spaces. \ep

\rks The monopole space, mentioned earlier, can also be described as a
certain space of rational functions on the Riemann sphere. The {\em
Sen conjecture} (see \cite{sen}) makes certain predictions about the
$L^2$-cohomology of this space, in particular that there is a certain
$SL(2, \Z)$ action on it. Might this action be derivable from TQFT?

%%%%%%%%%%%%%%%%%%%%%%%%%%%%%%%%%%%%%%%%%%%%%%%%%%%%%%%%%%%%%%%%%%%%%

\bp Define RW invariants for singular hyperk\"ahler manifolds. \ep

\rks It might be helpful to be able to handle orbifolds, in order to
better understand the behaviour of RW invariants of hyperk\"ahler
quotients and quotients by finite groups. If the equivariant theory
goes through, might there be localisation formulae for the curvature
integrals?

%%%%%%%%%%%%%%%%%%%%%%%%%%%%%%%%%%%%%%%%%%%%%%%%%%%%%%%%%%%%%%%%%%%%%

\bp Compute the invariants for some non-K\"ahler holomorphic
symplectic manifolds. \ep

\rks It would be interesting to calculate RW invariants for
non-K\"ahler holomorphic symplectic manifolds, such as Douady spaces
of the Kodaira surface (in fact the invariants are zero for these) or
Guan's examples \cite{guan}, and to try to work out whether such
invariants are independent of choice (if any) of complex structure.

%%%%%%%%%%%%%%%%%%%%%%%%%%%%%%%%%%%%%%%%%%%%%%%%%%%%%%%%%%%%%%%%%%%%%

\bp Try to extend the invariants to other holonomy groups. \ep

\rks Why stop with quaternionic-K\"ahler manifolds? If we take the
viewpoint that Rozansky-Witten invariants are special characteristic
classes arising out of a reduction of structure group on a Riemannian
manifold, then it is natural to try to find analogues for all special
holonomy groups.

Here is the idea. First examine the symmetries of the curvature tensor
and decide what kind of graphical vertex it corresponds to. Decide on
what kind of graphs index the possible self-combinations of the
curvature and the other special structural forms (e.g. on a
$G_2$-manifold, the canonical $3$- and $4$-forms.) Then attempt to
compute the universal graphical relations amongst such combinations;
that is, find out what kind of graph homology is appropriate.

For example, on any K\"ahler manifold (holonomy $U(k)$) one is free to
use trivalent {\em trees} to combine the curvature with itself, and
the IHX relation is satisfied in cohomology. This is part of
Kapranov's $L_\infty$ structure. Graphs with one loop can be created
by taking trace, but to glue up more legs would require a bilinear
form, as in the holomorphic symplectic case. On a {\em Calabi-Yau}
manifold (holonomy $SU(k)$) there is a canonical volume form, which
would be pictured as a $k$-valent vertex. Whether or not the
appropriate graph complex is ever rich enough to be useful, we do not know.

%%%%%%%%%%%%%%%%%%%%%%%%%%%%%%%%%%%%%%%%%%%%%%%%%%%%%%%%%%%%%%%%%%%%%

\bp Prove that $S^6$ is not complex. \ep

\rks It is a famous problem to show that the $6$-sphere, though
almost-complex, has no (integrable) complex structure. Perhaps the
$L_\infty$ structure on the Dolbeault complex of a complex manifold
(due to Kapranov, and a basic ingredient in RW theory) might provide a
new approach to a contradiction.

%%%%%%%%%%%%%%%%%%%%%%%%%%%%%%%%%%%%%%%%%%%%%%%%%%%%%%%%%%%%%%%%%%%%%

\bp Explore Kontsevich's generalisation of the theory. \ep

\rks Kontsevich gave a tremendous generalisation of the RW formalism,
taking in foliated {\em real} manifolds with a transverse symplectic
form, and in particular flat symplectic fibre bundles. This
construction is completely unexplored; perhaps a sheaf-theoretic
approach could help in calculating some examples, starting with the
case of surface bundles over surfaces.

%%%%%%%%%%%%%%%%%%%%%%%%%%%%%%%%%%%%%%%%%%%%%%%%%%%%%%%%%%%%%%%%%%%%%

\bp Investigate the use of higher graph cohomology. \ep

\rks Kapranov's $L_\infty$ structure allows the construction of weight
systems in higher graph cohomology, that is, numerical evaluations of
graphs with vertices of valence bigger than three. Little is known
about the graph homology spaces in positive degree, primarily because
the usual Lie algebraic methods are useless for studying them.

%%%%%%%%%%%%%%%%%%%%%%%%%%%%%%%%%%%%%%%%%%%%%%%%%%%%%%%%%%%%%%%%%%%%%

\bp Investigate the functoriality of RW invariants. \ep

\rks In contrast to the usual theory of characteristic classes, the RW
theory does not have a theory of {\em classifying spaces}, unless one
uses Kontsevich's Gelfand-Fuchs cohomology framework. A basic problem
is that hyperk\"ahler manifolds don't form a reasonable category
(there are almost no morphisms), and so comparing the invariants of
different manifolds is very hard. Is there any way to get around this?

%%%%%%%%%%%%%%%%%%%%%%%%%%%%%%%%%%%%%%%%%%%%%%%%%%%%%%%%%%%%%%%%%%%%%

\bp Compare Chern-Simons theory and RW theory as sigma-models. \ep

\rks Chern-Simons theory is a TQFT defined by Witten \cite{wittencs}
whose partition function is given by an integral over the space of all
connections on the trivial principal $G$-bundle over a
$3$-manifold. Such connections can be thought of as pullbacks of the
universal connection on the bundle $EG \rightarrow BG$ via smooth maps
$M \rightarrow BG$, and so in these terms CS theory resembles the RW
theory with $X=BG$ (compare this with the remarks in section 4 about
the loop space). The CS integrand can be defined by extending the map
$M \rightarrow BG$ over a $4$-manifold $W$ whose boundary is $M$, and
computing the integral over $W$ of the pullback of the fundamental
(Pontrjagin) $4$-form on $BG$. In the RW case one could try to do
something similar, at least for a null-homologous $M \rightarrow X$,
by using the fundamental $4$-form $\omega \bar\omega$. Can these
analogies be made more precise, and if so is there a fruitful exchange
of techniques?

%%%%%%%%%%%%%%%%%%%%%%%%%%%%%%%%%%%%%%%%%%%%%%%%%%%%%%%%%%%%%%%%%%%%%

\bp What are the trajectories of the RW action functional? \ep

\rks In Chern-Simons theory, the action functional can be viewed as a
Morse function on the space of $G$-connections on the $3$-manifold
$M$. Its critical points are the flat connections, and the gradient
flow lines (corresponding to instantons on $M \times \R$) can be used
to define the Floer homology of $M$. In RW theory, the action
functional is defined on smooth maps $\phi: M^3 \rightarrow X$, and
its stationary points are the constant maps. Can one make sense of the
trajectories of the functional, thought of as maps $: M \times \R
\rightarrow X$, and does this lead to a Floer homology of some kind?

%%%%%%%%%%%%%%%%%%%%%%%%%%%%%%%%%%%%%%%%%%%%%%%%%%%%%%%%%%%%%%%%%%%%%

\bp Can RW weight systems detect orientation? \ep

\rks There are many problems on the subject of the RW TQFT, and its
interaction with the theories of Vassiliev and quantum
invariants. These will be addressed elsewhere \cite{rsw01}, but
perhaps the most straightforward is as follows. The Vassiliev algebra
${\mathcal A}$ has an involution induced by reversing orientations of
the preferred circles of all spanning diagrams. An open question is
whether this involution is the identity or not; equivalently, whether
Vassiliev invariants cannot (or can) detect the orientation of
knots. Weight systems arising from Lie algebras are incapable of
distinguishing between diagrams and their reverses, but it is not at
all clear how to prove that RW weight systems share this
property. Thus it is conceivable that they could be used to find
non-reversible diagrams, if indeed such things exist.

\medskip

{\bf Acknowledgements}\qua The first author was supported by NSF Grant
DMS-0103922 and JSPS fellowship S-01037; the second by New College,
Oxford and by RIMS. We thank the organisers very much for our
invitations to RIMS and for their hard work in setting up such a
stimulating and enjoyable conference.

%%%%%%%%%%%%%%%%%%%%   End of main body of article
%
%                             References
%

\Addresses\recd

\end{document}

%% file: gtmonout.tex
\def\ifplaintex{\expandafter\ifx\csname documentclass\endcsname\relax}

\ifplaintex 
\else
\fi

\def\gtm{{\mathsurround=0pt\it $\cal G\mskip-2mu$eometry \&\ 
$\cal T\!\!$opology $\cal M\mskip-1mu$onographs}}    %  for monographs

\def\gtp{{\mathsurround=0pt\it $\cal G\mskip-2mu$eometry \&\ 
$\cal T\!\!$opology $\cal P\!$ublications}}  % GT publications

\def\recd{{\small Received:\qua\receiveddate\ifx\reviseddate\relax
\else\qquad Revised:\qua\reviseddate\fi\par}} 

\def\volumenumber#1{\def\thevolumenumber{#1}}
\def\volumeyear#1{\def\thevolumeyear{#1}}
\def\volumename#1{\def\thevolumename{#1}}
\def\papernumber#1{\def\thepapernumber{#1}}
\def\pagenumbers#1#2{\def\startpage{#1}\def\finishpage{#2}}
\def\published#1{\def\publishdate{#1}}
\def\received#1{\def\receiveddate{#1}}
\def\revised#1{\def\reviseddate{#1}}
\def\accepted#1{\def\accepteddate{#1}}

\def\asciiaddress#1{\def\theasciiaddress{#1}}

\long\def\asciiabstract#1{\long\def\theasciiabstract{#1}}
\def\asciikeywords#1{\def\theasciikeywords{#1}}

\let\\\par
\let\thevolumenumber\relax\let\thepapernumber\relax
\let\thevolumeyear\relax\let\startpage\relax
\let\finishpage\relax\let\publishdate\relax\let\receiveddate\relax
\let\reviseddate\relax\let\accepteddate\relax\let\theasciititle\relax
\let\theasciiauthors\relax\let\theasciiaddress\relax
\let\theasciiabstract\relax\let\theasciikeywords\relax

\let\theerratum\relax\let\theasciiemail\relax
\let\theshortauthors\relax\let\theshorttitle\relax

\def\startpage{1}\def\finishpage{15}\def\thepapernumber{77}

\volumenumber{2}
\volumename{Proceedings of the Kirbyfest}
\volumeyear{1999}

\long\def\maketitlep{   % start of definition of \maketitlep

\count0=\startpage

\gtm\nl        %   GT mongraphs (top left) 
{\small Volume \thevolumenumber: \thevolumename\nl 
\ifx\theerratum\relax\else Erratum \erratumnumber\nl\fi
Pages \startpage--\finishpage\nl}

\vglue 0.1truein   % top margin

{\parskip=0pt\leftskip 0pt plus 1fil\def\\{\par\smallskip}{\ifplaintex\large
\else\Large\fi\bf\thetitle}\par\medskip}   
\vglue 0.05truein 

{\parskip=0pt\leftskip 0pt plus 1fil\def\\{\par}{\sc\theauthors}
\par\medskip}%
 
\vglue 0.03truein 

{\small\leftskip 25pt\rightskip 25pt{\bf Abstract}\stdspace\theabstract

{\bf AMS Classification}\stdspace\theprimaryclass
\ifx\thesecondaryclass\relax\else; \thesecondaryclass\fi\par
{\bf Keywords}\stdspace \thekeywords\par}\vglue 7pt

}   % end of definition of \maketitlep

\font\phead=cmsl9 scaled 950
\font=cmsl9 scaled 1050
\font\pnum=cmbx10 scaled 913
\font\lnum=cmbx10 
\font\pfoot=cmsl9 scaled 950
\font=cmsl9 scaled 1050
\ifplaintex
\headline{\vbox to 0pt{\vskip -4.5mm\line{\small\phead\ifnum
\count0=\startpage ISSN 1464-8997 (on line)
1464-8989 (printed) \hfill {\pnum\folio}\else\ifodd\count0\def\\{ }% 
\ifx\theshorttitle\relax\thetitle\else\theshorttitle\fi\hfill{\pnum\folio}
\else\def\\{ and }{\pnum\folio}\hfill\ifx\theshortauthors\relax\theauthors
\else\theshortauthors\fi\fi\fi}\vss}}
\footline{\vbox to 0pt{\vglue 0mm\line{\small\pfoot\ifnum\count0=\startpage
Published \publishdate:\qua\copyright\ \gtp\hfill\else
\gtm, Volume \thevolumenumber\ (\thevolumeyear)\hfill\fi}\vss
}}
\else
\makeatletter
\def\@oddhead{{\small\lhead\ifnum\count0=\startpage ISSN 1464-8997 (on line)
1464-8989 (printed) \hfill {\lnum\number\count0}\else\ifodd\count0
\def\\{ }\ifx\theshorttitle\relax \thetitle \else\theshorttitle\fi\hfill
{\lnum\number\count0}\else\def\\{ and }{\lnum\number\count0}
\hfill\ifx\theshortauthors\relax 
\theauthors\else\theshortauthors\fi\fi\fi}}\def\@evenhead{@oddhead}
\def\@oddfoot{\small\lfoot\ifnum\count0=\startpage Published \publishdate:\qua\copyright\ \gtp\hfill\else
\gtm, Volume \thevolumenumber\ (\thevolumeyear)\hfill\fi}
\def\@evenfoot{@oddfoot}
\makeatother
\fi

\let\maketitlepage\maketitlep
\let\makeshorttitle\maketitlepage
\let\maketitle\maketitlepage

\newwrite\gtoutfile
\long\gdef\makeheadfile{  %%% start of definition of \makeheadfile
{\def\\{, }\def\s{ }
\immediate\openout\gtoutfile head.xxx
\immediate\write\gtoutfile{To: math@arxiv.org}
\immediate\write\gtoutfile{Subject: put OR rep NNNNN:ppppp}
\immediate\write\gtoutfile{--text follows this line--}
\immediate\write\gtoutfile{Proxy-for: \ifx\theasciiauthors\relax
\theauthors\else\theasciiauthors\fi\s<\ifx\theasciiemail\relax\theemail\else\theasciiemail\fi>}
\immediate\write\gtoutfile{\noexpand\\}
\immediate\write\gtoutfile{Authors: \ifx\theasciiauthors\relax
\theauthors\else\theasciiauthors\fi}
{\def\\{ }\immediate\write\gtoutfile{Title: \ifx\theasciititle\relax
\thetitle\else\theasciititle\fi}}
\immediate\write\gtoutfile{Subj-class: GT or SG, GR etc}
\immediate\write\gtoutfile{MSC-class: \theprimaryclass\ifx\thesecondaryclass\relax\else, \thesecondaryclass\fi}
\immediate\write\gtoutfile{Journal-ref: Geom. Topol. Monogr. \thevolumenumber\s
(\thevolumeyear) \startpage-\finishpage}
\immediate\write\gtoutfile{Comments: Published by Geometry and Topology Monographs at}
\immediate\write\gtoutfile{\s\s\s  http://www.maths.warwick.ac.uk/gt/GTMon\thevolumenumber/paper\thepapernumber.abs.html}
\immediate\write\gtoutfile{\noexpand\\}
\immediate\write\gtoutfile{}
\ifx\theasciiabstract\relax
\immediate\write\gtoutfile{\theabstract}\else
\immediate\write\gtoutfile{\theasciiabstract}\fi
\immediate\write\gtoutfile{}
\immediate\write\gtoutfile{\noexpand\\}
\immediate\write\gtoutfile{}
\immediate\closeout\gtoutfile}}  %%% end of definition of \makeheadfile

\def\maketitlepage{\maketitlep\makeheadfile}
\let\makeshorttitle\maketitlepage
\let\maketitle\maketitlepage